\documentclass[a4paper,12pt]{article}

\usepackage[left=2cm,right=2cm, top=2cm,bottom=2cm,bindingoffset=0cm]{geometry}

\usepackage{verbatim}
\usepackage{amsmath}
\usepackage{amsthm}
\usepackage{amssymb}
\usepackage{delarray}
\usepackage{mathrsfs}
\usepackage{graphicx}

\theoremstyle{plain}

\theoremstyle{definition}

\theoremstyle{remark}

\DeclareMathOperator{\rank}{rank}
 \allowdisplaybreaks
\usepackage{titlesec}
\titleformat{\section}
{\normalfont\bfseries\centering}{}{0em}{}

\numberwithin{equation}{section}
\begin{document}

\begin{center}
{\bf On recovering Sturm-Liouville operators with frozen argument from the spectrum}
\end{center}
\begin{center}
{\bf N.\,P.~Bondarenko\footnote{Department of Applied Mathematics, Samara National Research University; Department of Mathematics, Saratov
National Research State University, email: bondarenkonp@info.sgu.ru}, S.\,A.~Buterin\footnote{Department of Mathematics, Saratov National
Research State University, email: buterinsa@info.sgu.ru} and S.\,V.~Vasiliev\footnote{Department of Mathematics, Saratov National Research
State University, email: altrair@mail.ru}}
\end{center}

\begin{center}
\parbox[c]{130mm}{\small {\bf Abstract.} We consider second order linear differential operators possessing a term depending on the unknown
function with a fixed argument and study the uniqueness of recovering the operators from the spectrum. We also obtain a constructive
procedure for solving this inverse problem along with necessary and sufficient conditions of its solvability.
\\


}
\end{center}

\section{1. Introduction}

Let $\{\lambda_n\}_{n\ge1}$ be the spectrum of the boundary value problem $L=L(q(x),a,\alpha,\beta)$ of the form
\begin{equation}\label{1.1}
\ell y:=-y''(x)+q(x)y(a)=\lambda y(x), \quad 0<x<\pi,
\end{equation}
\begin{equation}\label{1.2}
y^{(\alpha)}(0)=y^{(\beta)}(\pi)=0,
\end{equation}
where $\lambda$ is the spectral parameter, $q(x)$ is a complex-valued function in $L_2(0,\pi)$ and $\alpha,\beta\in\{0,1\}.$ Let also
$k:=\pi/a\in{\mathbb N}.$ The case $k\notin{\mathbb N}$ requires a separate investigation. We call $\ell$ the Sturm-Liouville operator with
{\it frozen argument}.

In this paper we study an inverse spectral problem for $L.$ Inverse problems of spectral analysis consist in recovering operators from
their spectral characteristics. The greatest success in the inverse spectral theory has been achieved for the classical Sturm-Liouville
operator (see, e.g., \cite{Bor}--\cite{RS}) and afterwards for higher order differential operators \cite{Beals}--\cite{Yur-2}. For example,
it is known \cite{Bor} that the potential $q(x)$ is uniquely determined by the spectra of two boundary value problems for one and the same
classical Sturm-Liouville equation
$$
-y''(x)+q(x)y(x)=\lambda y(x), \quad 0<x<\pi,
$$
with one common boundary condition. For differential operators with frozen argument as well as for other classes of non-local operators the
classical methods of inverse spectral theory do not work and there are only few results in this direction, which do not form a
comprehensive picture. Some aspects of inverse spectral theory for differential operators with frozen argument were studied in
\cite{AlbHryNizh}--\cite{BV-3}. For example, the authors of \cite{AlbHryNizh} studied the case of a real-valued $q(x)$ and $a=\pi$ with the
Dirichlet boundary condition in the point $x=0$ and the nonlocal condition $y'(\pi)+\int_0^\pi y(x)q(x)\,dx=0$ depending on the potential.
In \cite{BV-1}--\cite{BV-3} an inverse problem was studied for $L$ with some particular values of $\alpha$ and $\beta.$

Differential operators with frozen argument can be classified as a special case of differential operators with deviating argument
\cite{Nor}--\cite{Vla}, which have many application in natural sciences and engineering. In particular, inverse spectral problems for
differential operators with constant delay were studied in \cite{Pik}--\cite{BY}, while the case of integral delay was investigated in
\cite{Yur-3}--\cite{BS} and other papers. In the present paper we study the following inverse problem.

\medskip
{\bf Inverse Problem 1.1.} Given $\{\lambda_n\}_{n\ge1},$ $a,$ $\alpha$ and $\beta,$ find $q(x).$

\medskip
We establish, in particular, that for certain values of $\alpha,$ $\beta$ and $k,$ unlike to the classical case, the specification of the
spectrum is sufficient for unique determination of the potential. We refer to this case when the uniqueness holds as non-degenerate one,
while to the opposite case we refer as degenerate one. In the degenerate case we describe classes of iso-spectral potentials and suggest an
additional restriction on the potential under which the uniqueness holds. For example, for $k>1$ the uniqueness in the degenerate case can
be achieved by restricting the class of potentials by the additional assumption
\begin{equation}\label{1.3}
q(a-t)=K(q(a+t)), \quad 0<t<a,
\end{equation}
where $K$ is some operator in $L_2(0,a)$ such that $I+K$ is invertible and $I$ is the identity operator. In
particular, if $K\equiv I,$ then the condition \eqref{1.3} is equivalent to the evenness of the potential
with respect to the point $a,$ i.e. $q(a-t)=q(a+t)$, $0<t<a.$ We note that the case of odd potentials with
respect to the point $a$ is not covered by the condition \eqref{1.3} and not eligible. Indeed, as can be seen
in the proof of Theorem~2.1 below, if $K=-I,$ then the spectrum of $L(q(x),\pi/2,0,0)$ coincides with the
spectrum of $L(0,a,0,0)$ and, hence, carries no information on $q(x).$ The case of constant $K$ (i.e. when
$K(f)$ does not depend on $f)$ corresponds to a priori specification of $q(x)$ on the subinterval $(0,a).$
There may be used also other than \eqref{1.3} restrictions guarantying the uniqueness of solution of Inverse
Problem~1.1 in the degenerate case, which, generally speaking, depend on the parameters $\alpha,$ $\beta$ and
$k.$

We also obtain a characterization of the spectrum of the boundary value problem \eqref{1.1}, \eqref{1.2}. In
other words, we obtain conditions on $\{\lambda_n\}_{n\ge1}$ that are necessary and sufficient for the
solvability of Inverse Problem~1.1. The related proof is constructive and gives algorithms for solving the
inverse problem. Thus, we, actually, obtain a complete description of all possible situations in the case of
natural $k.$ Using our approach one can also investigate the case of rational $k,$ which promises to have
more complicated description of degenerate and non-degenerate subcases (see Remark~4.2 below).

The paper is organized as follows. In the next section we reduce our inverse problem to the so-called main equation and investigate its
solvability. With accordance to this, we introduce and study the degenerate and non-degenerate cases. In Section~3 we establish properties
of the spectrum for various combinations of values of the parameters $\alpha,$ $\beta$ and $k.$ In Section 4 we prove the uniqueness
theorem and provide constructive procedures for solving the inverse problem along with necessary and sufficient conditions of its
solvability both in the degenerate and non-degenerate cases. Moreover, the set of all iso-spectral potentials in the degenerate case is
described.

\section{2. Main equation of the inverse problem. Degenerate and non-degenerate cases}

Let $C(x,\lambda), \ S(x,\lambda)$ be solutions of equation (\ref{1.1}) under the initial conditions
$$
C(a,\lambda)=S'(a,\lambda)=1, \quad S(a,\lambda)=C'(a,\lambda)=0.
$$
It is easy to check that
\begin{equation}\label{2.1}
C(x,\lambda)=\cos\rho(x-a)+\int\limits_a^x\frac{\sin\rho(x-t)}{\rho}q(t)\,dt,
\end{equation}
\begin{equation}\label{2.2}
S(x,\lambda)=\frac{\sin\rho(x-a)}{\rho}.
\end{equation}
Clearly, eigenvalues of $L$ coincide with the zeros of its characteristics function
\begin{equation}\label{2.3}
\Delta_{\alpha,\beta}(\lambda)=
\begin{vmatrix}
C^{(\alpha)}(0,\lambda) & S^{(\alpha)}(0,\lambda) \\
C^{(\beta)}(\pi,\lambda) & S^{(\beta)}(\pi,\lambda)
\end{vmatrix}.
\end{equation}
Let $f(t)\in L_2(0,\pi).$ Recalling that $ka=\pi$ for a certain $k\in{\mathbb N},$ we introduce the following shift and involution
operators
\begin{equation}\label{2.3.1}
R_mf(t)=\left\{
\begin{array}{l}
\displaystyle\!\! f(t+(k-m)a) \;\; \text{for odd} \; m,\\[3mm]
\displaystyle\!\! f((k-m+1)a-t) \;\; \text{for even} \; m,
\end{array}\right.
\!\!Q_mf(t)=\left\{
\begin{array}{ll}
\displaystyle\!\! f(t+(m-1)a) \;\; \text{for odd} \; m,\\[3mm]
\displaystyle\!\! f(ma-t) \;\; \text{for even} \; m,
\end{array}\right.
\end{equation}
where $t\in(0,a)$ and $m=\overline{1,k}.$ Consider the operators $R,Q:L_2(0,\pi)\to(L_2(0,a))^k$ determined by the formulae
\begin{equation} \label{2.4}
Rf:=(R_1f,R_2f,\ldots,R_kf)^T, \quad Qf:=(f,Q_2f,\ldots,Q_kf)^T,
\end{equation}
where $T$ is the transposition sign. Obviously, the operators $R$ and $Q$ are invertible.

\medskip
{\bf Theorem 2.1. }{\it The characteristic function $\Delta_{\alpha,\beta}(\lambda)$ of the problem $L$ has
the form
\begin{equation} \label{2.5}
\Delta_{\alpha,\alpha}(\lambda)=\rho^{2\alpha}\Big(\frac{\sin\rho\pi}{\rho}+
\int\limits_0^{\pi}W_{\alpha,\alpha}(t)\frac{\cos\rho t}{\rho^2}\,dt\Big),\quad W_{\alpha,\alpha}(t)\in
L_2(0,\pi), \; \int\limits_0^{\pi}W_{0,0}(t)\,dt=0,
\end{equation}
if $\alpha=\beta$ and
\begin{equation}\label{2.6}
\Delta_{\alpha,\beta}(\lambda)=(-1)^\alpha\cos\rho\pi+\int\limits_0^{\pi}W_{\alpha,\beta}(t)\frac{\sin\rho
t}{\rho}dt, \quad W_{\alpha,\beta}(t)\in L_2(0,\pi),
\end{equation}
if $\alpha\ne\beta.$ Moreover, the function $W_{\alpha,\beta}(t)$ has the form
\begin{equation}\label{2.7}
W_{\alpha,\beta}(t) = \frac{(-1)^{\alpha\beta}}2 Q^{-1}A_{\alpha,\beta} Rq(t),
\end{equation}
where $A_{\alpha,\beta}$ is a square three-diagonal matrix of order $k$ having the form
\begin{equation}\label{2.8}
A_{\alpha,\beta}=\begin{pmatrix}
1 & b & 0 & 0 & 0 & \cdots & 0\\
c & 0 & b & 0 & 0 & \cdots & 0\\
0 & c & 0 & b & 0 & \cdots & 0\\
\vdots & \vdots & \ddots & \ddots & \ddots & \ddots & \vdots\\
0 & 0 & \cdots & c & 0 & b & 0\\
0 & 0 & 0 & \cdots & c & 0 & b\\
0 & 0 & 0 & 0 & \cdots & c & c\\
\end{pmatrix}
\end{equation}
for $k>1$ and
\begin{equation}\label{2.8-1}
A_{\alpha,\beta}=2(-1)^{\alpha(\beta+1)}\delta_{1,\beta}
\end{equation} for $k=1,$ where
$\delta_{1,\beta}$ is the Kronecker delta.

In \eqref{2.8} each subdiagonal consists of equal elements and all the elements of the main diagonal, except
the first and last ones, vanish. Moreover,
\begin{equation} \label{2.9}
b=(-1)^{\alpha+\beta}, \quad c=(-1)^{1+\beta}.
\end{equation}
}

\medskip
{\it Proof.} Let for definiteness $k>1.$ Consider first the situation when $\alpha=\beta=0$. Substituting
\eqref{2.1} and \eqref{2.2} into \eqref{2.3} and using the formula
$$
2\sin\rho t\sin\rho(\pi-a)=\cos\rho(\pi-a-t)-\cos\rho(\pi-a+t),
$$
we obtain the representation
$$
\Delta_{0,0}(\lambda)=\frac{\sin\rho\pi}{\rho}
+\frac{1}{2\rho^2}\int\limits_0^{a}q(t)\Big(\cos\rho(\pi-a-t)-\cos\rho(\pi-a+t)\Big)\,dt+
$$
$$
+\frac{1}{2\rho^2}\int\limits_{a}^\pi q(t)\Big(\cos\rho(\pi-a-t)-\cos\rho(\pi+a-t)\Big)\,dt.
$$
Changing the variables of integration, we get
$$
\Delta_{0,0}(\lambda)=\frac{\sin\rho\pi}{\rho}
+\frac{1}{2\rho^2}\int\limits_{\pi-2a}^{\pi-a}q(\pi-a-t)\cos\rho t\,dt-
\frac{1}{2\rho^2}\int\limits_{\pi-a}^{\pi}q(t-\pi+a)\cos\rho t\,dt+
$$
$$
+\frac{1}{2\rho^2}\int\limits_{-a}^{\pi-2a} q(\pi-a-t)\cos\rho t\,dt -\frac{1}{2\rho^2}\int\limits_{a}^\pi q(\pi+a-t)\cos\rho t\,dt.
$$
Thus, we  arrive at \eqref{2.5} for $\alpha=\beta=0,$ where the function $W_{0,0}(x)$ has the form
\begin{equation}\label{dd:1}
\displaystyle W_{0,0}(t)=\frac{1}{2} \left\{
\begin{array}{l}
\displaystyle q(\pi-a+t)+q(\pi-a-t), \quad t\in(0,a),\\[3mm]
\displaystyle -q(\pi+a-t)+q(\pi-a-t), \quad t\in(a,\pi-a),\\[3mm]
\displaystyle -q(\pi+a-t)-q(t-\pi+a), \quad t\in(\pi-a,\pi).
\end{array}
\right.
\end{equation}
We note that the last equality in \eqref{2.5} follows from the entireness of the function
$\Delta_{0,0}(\lambda).$

Proceeding analogously for the other combinations of $\alpha,\beta\in\{0,1\}$ and taking \eqref{2.9} and $\pi=ka$ into account one can
obtain the following general representation
\begin{equation}\label{allW}
W_{\alpha,\beta}(t)=\frac{(-1)^{\alpha\beta}}{2} \left\{
\begin{array}{l}
\displaystyle q((k-1)a+t)+bq((k-1)a-t), \quad t\in(0,a),\\[3mm]
\displaystyle cq((k+1)a-t)+bq((k-1)a-t), \quad t\in(a,(k-1)a),\\[3mm]
\displaystyle c\Big(q((k+1)a-t)+q(t-(k-1)a)\Big), \quad t\in\Big((k-1)a,ka\Big).
\end{array}
\right.
\end{equation}
Acting by the operator $Q$ on the both sides of \eqref{allW} and using \eqref{2.4} we get
\begin{equation}\label{sys}
QW_{\alpha,\beta}(t)=\frac{(-1)^{\alpha\beta}}{2} \left(
\begin{array}{c}
R_1q(t)+bR_2q(t)\\[2mm]
cR_1q(t)+bR_3q(t)\\[2mm]
cR_2q(t)+bR_4q(t)\\[2mm]
cR_3q(t)+bR_5q(t)\\[2mm]
\ldots\\[2mm]
cR_{k-2}q(t)+bR_kq(t)\\[2mm]
c(R_{k-1}q(t)+R_kq(t))
\end{array}\right).
\end{equation}
Taking into account \eqref{2.8} and invertibility of $Q$ we arrive \eqref{2.7}.  The case $k=1$ is easier and
can be treated in a similar way. $\hfill\Box$

\medskip
Equation \eqref{2.7} is called the {\it main equation} of Inverse Problem 1.1. The following lemma gives a formula for
$\det{A_{\alpha,\beta}}.$

\medskip
{\bf Lemma 2.1. }{\it Let $k>1.$ The determinant of the matrix $A_{\alpha,\beta}$ of the form \eqref{2.8} can
be calculated by the formula
\begin{equation}\label{charfun:9}
\det A_{\alpha,\beta} = \left\{
\begin{array}{l}
\displaystyle (-bc)^{(k-1)/2}(1+c), \quad \text{if k is odd},\\[3mm]
\displaystyle (-b)^{k/2-1}c^{k/2}(1-b), \quad \text{if k is even}.
\end{array}\right.
\end{equation}
}

\bigskip
{\it Proof.} In $\det A_{\alpha,\beta}$ sequentially for $j=\overline{3,k}$ subtracting $(j-2)$-th row
multiplied with $c/b$ from the $j$-th row and then by the analogous way zeroizing the first column except for
its last element, we obtain
$$
\det A_{\alpha,\beta}=\left|\begin{array}{ccccccc}
1 & b & 0 & 0 & 0 & \cdots & 0\\
c & 0 & b & 0 & 0 & \cdots & 0\\
p_3 & 0 & 0 & b & 0 & \cdots & 0\\
\vdots & \vdots & \ddots & \ddots & \ddots & \ddots & \vdots\\
p_{k-2} & 0 & \cdots & 0 & 0 & b & 0\\
p_{k-1} & 0 & 0 & \cdots & 0 & 0 & b\\
p_k & 0 & 0 & 0 & \cdots & 0 & c\\
\end{array}\right|=
\left|\begin{array}{ccccccc}
0 & b & 0 & 0 & 0 & \cdots & 0\\
0 & 0 & b & 0 & 0 & \cdots & 0\\
0 & 0 & 0 & b & 0 & \cdots & 0\\
\vdots & \vdots & \ddots & \ddots & \ddots & \ddots & \vdots\\
0 & 0 & \cdots & 0 & 0 & b & 0\\
0 & 0 & 0 & \cdots & 0 & 0 & b\\
p & 0 & 0 & 0 & \cdots & 0 & c\\
\end{array}\right|,
$$
where $p_{2j+1}=(-c/b)^j,$ $p_{2j+2}=c(-c/b)^j,$ $j\ge1,$ and $p=p_k-cp_{k-1}/b.$ Thus, we arrive at
$$
\det A_{\alpha,\beta}=(-b)^{k-1}\Big(p_k-\frac{c}{b}p_{k-1}\Big)=(-b)^{k-1}\left\{
\begin{array}{l}
\displaystyle \Big(-\frac{c}{b}\Big)^{(k-1)/2}(1+c), \quad \text{if k is odd},\\[3mm]
\displaystyle c\Big(-\frac{c}{b}\Big)^{k/2-1}\Big(1-\frac1b\Big), \quad \text{if k is even},
\end{array}\right.
$$
which is equivalent to \eqref{charfun:9}. $\hfill\Box$

\medskip
We shall reduce Inverse Problem~1.1 to solving the main equation \eqref{2.7}. While the case $k=1$ is
trivial, Lemma~2.1 gives us the knowledge about when equation \eqref{2.7} is uniquely solvable with respect
to $q(x)$ when $k>1.$ In accordance with the existence of two possibilities, we highlight two cases: {\it
degenerate} and {\it non-degenerate} ones, depending on whether $\det A_{\alpha,\beta}=0$ or $\det
A_{\alpha,\beta}\ne0,$ respectively. Since, by virtue of \eqref{2.9}, we have
$$
1+c=1+(-1)^{1+\beta}, \quad  1-b=1-(-1)^{\alpha+\beta},
$$
according to Lemma~2.1 and \eqref{2.8-1} the degenerate case occurs when one of the following groups of
conditions is fulfilled:
\begin{equation}\label{deg}
\left.\begin{array}{rl}
\text{(i)}& \alpha=\beta=0;\\[3mm]
\text{(ii)}& \alpha=1,\ \beta=0 \text{ and } k \text{ is odd;}\\[3mm]
\text{(iii)}& \alpha=\beta=1 \text{ and } k \text{ is even;}\\[3mm]
\end{array}\right\}
\end{equation}
while the non-degenerate case includes the remaining groups of conditions:
\begin{equation}\label{non-deg}
\left.\begin{array}{rl}
\text{(iv)}& \alpha=0, \;\beta=1;\\[3mm]
\text{(v)}& \alpha=1,\ \beta=0 \text{ and } k \text{ is even;}\\[3mm]
\text{(vi)}& \alpha=\beta=1 \text{ and } k \text{ is odd.}\\[3mm]
\end{array}\right\}
\end{equation}

We note that each time when solving Inverse Problem~1.1 for $L=L(q(x),\pi/k,\alpha,\beta)$ we shall assume that $k\in{\mathbb N}$ along
with $\alpha,\beta\in\{0,1\}$ are fixed and known a priori. If the corresponding problem $L$ belongs to the non-degenerate case we shall
prove in Section~4 the uniqueness of solution and obtain a constructive procedure for solving the inverse problem along with necessary and
sufficient conditions of its solvability. The latter are, actually, equivalent to the full characterization of the spectrum of $L$ in terms
of asymptotics. For the degenerate case we shall also obtain the characterization of the spectrum, which besides the asymptotics will
include an additional degeneration condition. In the class of potentials satisfying the restriction \eqref{1.3} we shall prove the
uniqueness of solution and obtain a constructive procedure for solving the inverse problem also in the degenerate case.

\section{3. Properties of the spectrum}

As was mentioned in the end of the preceding section, in the degenerate case the characterization of the spectrum has to include a certain
degeneration condition, which is, in turn, connected with some structural property of the function $W_{\alpha,\beta}(x).$ The following
lemma gives such structural properties for all groups of conditions representing the degenerate case.

\medskip
{\bf Lemma 3.1. }{\it In the degenerate case the function $W_{\alpha,\beta}(t)$, determined by \eqref{allW}, satisfies one of the following
equalities:

(i) For $\alpha=\beta=0:$
\begin{equation}\label{dd:3}
\sum_{j = 0}^{[(k-1)/2]} W_{0,0}( 2ja + t) + \sum_{j = 1}^{[k/2]} W_{0,0} ( 2ja - t )  = 0 \quad \text{a.e. on}\; (0, a);
\end{equation}

(ii) For $\alpha=1,$ $\beta=0$ and odd $k:$
\begin{equation}\label{nd:3}
\sum_{j = 0}^{(k-1)/2} (-1)^j W_{1,0}( 2ja  + t ) =\sum_{j = 1}^{(k-1)/2} (-1)^j W_{1,0} ( 2ja  - t) \quad \text{a.e. on}\; (0, a);
\end{equation}

(iii) For $\alpha=\beta=1$ and even $k:$
\begin{equation}\label{nn:3}
\sum_{j = 0}^{k/2-1} (-1)^jW_{1,1}( 2ja  + t ) + \sum_{j = 1}^{k/2} (-1)^jW_{1,1}(2ja - t )  = 0, \quad \text{a.e. on}\; (0, a).
\end{equation}
Here $[x]$ denotes the entire part of $x.$}

\medskip
{\it Proof.} Let us prove \eqref{dd:3}. For $k=1$ it is obvious, let $k>1.$ By virtue of \eqref{dd:1}, for $1\le s\le[k/2]-1$ and $t\in
(0,a)$ we have
$$
W_{0,0}(t) + \sum_{j = 1}^{s} \Big( W_{0,0}(2ja+t)+W_{0,0}(2ja-t)\Big) =\frac12\Big(q((k-2s-1)a+t)+q((k-2s-1)a-t)\Big).
$$
Using this formula for $s=[k/2]-1$ along with the relation
$$
W_{0,0}((k-1)a+t)+W_{0,0}((k-1)a-t)=-\frac12\Big(q(2a+t)+q(2a-t)\Big), \quad t\in (0,a),
$$
for odd $k$ and the relation
$$
W_{0,0}(ka-t)=-\frac12\Big(q(a-t)+q(a+t)\Big), \quad t\in(0,a),
$$
for even $k$, we obtain \eqref{dd:3}. Analogously, using \eqref{allW} one can prove \eqref{nd:3} and \eqref{nn:3}. $\hfill\Box$

\medskip
From the proof of Lemma 2.1 it follows that in the degenerate case $\rank{A_{\alpha,\beta}}=k-1.$ This means that obeying only the
restriction imposed by Lemma~3.1 the function $W_{\alpha,\beta}(t)\in L_2(0,\pi)$ can be arbitrary in the rest. Actually, we prove this
later when studying the inverse problem. We note that the last condition in \eqref{2.5} also follows from \eqref{dd:3}, which, inter alia,
will be obtained as a consequence of \eqref{dd:3} in the proof of Theorem~3.1 (see below).

In the non-degenerate case $\rank{A_{\alpha,\beta}}=k$ and similarly one can prove the following lemma.

\medskip
{\bf Lemma 3.2. }{\it In the non-degenerate case the following relations hold:

(iv) For $\alpha=0$ and $\beta=1:$
$$
q(x)=W_{0,1}(x)+\sum_{j=1}^{(k-1)/2}\Big(W_{0,1}(2ja+x)-W_{0,1}(2ja-x)\Big), \quad x\in(0,a),
$$
if $k$ is odd, and
$$
q(x)=\sum_{j=1}^{k/2}\Big(W_{0,1}((2j-1)a+x)-W_{0,1}((2j-1)a-x)\Big), \quad x\in(0,a),
$$
if $k$ is even;

(v) For $\alpha=1,$ $\beta=0$ and even $k:$
$$
q(x)=\sum_{j=1}^{k/2}(-1)^j\Big(W_{1,0}((k+1-2j)a-x)+W_{1,0}((k+1-2j)a+x)\Big), \quad x\in(0,a);
$$

(vi) For $\alpha=\beta=1$ and odd $k:$
$$
q(x)=(-1)^\frac{k+1}2\Big(W_{1,1}(x) + \sum_{j=1}^{(k-1)/2}(-1)^j\Big(W_{1,1}(2ja+x)+W_{1,1}(2ja-x)\Big)\Big), \quad x\in(0,a).
$$
}

\medskip
The following theorem describes the properties of the spectrum of the problem $L.$

\medskip
{\bf Theorem 3.1. }{\it The problem $L$ has a countable set of eigenvalues $\{\lambda_n\}_{n\ge1}$ of the form
\begin{equation}
\lambda_n=\Big(n-\frac{\alpha+\beta}{2}+\frac{\kappa_n}{n}\Big)^2,\quad \{\kappa_n\}\in l_2. \label{dd:5}
\end{equation}
Moreover, in the degenerate case, a part of the eigenvalues degenerates in the following sense:

(i) For $\alpha=\beta=0:$
\begin{equation}\label{dd:6}
\lambda_{kn}=(kn)^2, \quad n\in \mathbb N;
\end{equation}

(ii) For $\alpha=1,$ $\beta=0$ and odd $k:$
\begin{equation}\label{BBV:11}
\lambda_{k(n-1/2)+1/2}=k^2\Big(n-\frac12\Big)^2, \quad n \in \mathbb N;
\end{equation}

(iii) For $\alpha=\beta=1$ and even $k:$
\begin{equation}\label{BBV:9}
\lambda_{k(n-1/2)+1}=k^2\Big(n-\frac12\Big)^2, \quad n\in{\mathbb N}.
\end{equation}
}

\medskip
{Proof.} The existence of a countable set of eigenvalues of the form \eqref{dd:5} can be established by the
standard approach involving Rouch\'e's theorem (see, e.g., [3]). It remains to prove
\eqref{dd:6}--\eqref{BBV:9}. The case $k=1$ is trivial. Let for definiteness $\alpha=\beta=0$ and $k>1.$
Expand $W_{0,0}(t)$ into the Fourier series
\begin{equation}\label{dd:7}
W_{0,0}(t)=\sum_{n=0}^\infty a_n\cos nt, \quad a_n=\frac2\pi\int_0^\pi W_{0,0}(t) \cos n t \, dt.
\end{equation}
Substituting \eqref{dd:7} together with the relation
$$
W_{0,0}(2ja+t)+W_{0,0}(2ja-t)=2\sum_{n=0}^\infty a_n\cos2jna\cos nt, \quad j=\overline{1,[(k-1)/2]}, \quad
t\in(0,a),
$$
into \eqref{dd:3}, we obtain
\begin{equation}\label{dd:8}
\sum_{j=0}^{[(k-1)/2]} W_{0,0}(2ja+t)+\sum_{j=1}^{[k/2]} W_{0,0}(2ja-t)= \sum_{n = 0}^{\infty} a_n T_n \cos n t=0, \quad t\in(0,a),
\end{equation}
where
\begin{equation}\label{dd:8-1}
T_n=1+2\sum_{j=1}^{[(k-1)/2]}\cos2jna+\frac{1+(-1)^k}2\cos\pi n.
\end{equation}
Since $\cos2naj=\cos2na(k-j),$ we arrive at
\begin{equation} \label{dd:9}
T_n = \sum_{j=0}^{k-1}\cos2jna= \mbox{Re}\sum_{j = 0}^{k - 1} \exp(2ijna) =\left\{\begin{array}{l}k,\quad n/k+1\in{\mathbb N},\\[3mm]
0,\quad n/k+1\notin{\mathbb N},\end{array}\right.
\end{equation}
because for $n/k+1\notin{\mathbb N}$ we have $\exp(2ina)\ne1,$ while
$$
(1-\exp(2ina))\sum_{j = 0}^{k - 1} \exp(2ijna)=1-\exp(2ikna)=0.
$$
According to \eqref{dd:9} the relation \eqref{dd:8} takes the form
$$
\sum_{n=0}^{\infty} a_{kn}\cos knt=0, \quad t\in(0,a),
$$
By virtue of the minimality of the functional system $\{\cos knt\}_{n\ge0}$ in $L_2(0,\pi/k),$ we get $a_{nk}=0$ for $n+1\in \mathbb N.$
Using \eqref{2.5} and \eqref{dd:7} we obtain $\Delta_{0,0}((kn)^2)=0$ for $n\in \mathbb N$ and, hence, \eqref{dd:6} holds. The relations
\eqref{BBV:11} and \eqref{BBV:9} can be proven analogously.  $\hfill\Box$

\medskip
By the standard approach using Hadamard's factorization theorem (see, e.g., in [3]) one can prove the following assertion.

\medskip
{\bf Lemma 3.3. }{\it The specification of the spectrum $\{\lambda_n\}_{n\ge1}$ uniquely determines the characteristic function by the
formula
\begin{equation}\label{EQ}
\Delta_{\alpha,\beta}(\lambda) =\pi^{\delta_{\alpha,\beta}}(\lambda_1-\lambda)^{\alpha\beta}
\prod\limits_{n=1+\alpha\beta}^\infty\frac{\lambda_n-\lambda}{\Big(n-\frac{\alpha+\beta}2\Big)^2},
\end{equation}
where $\delta_{\alpha,\beta}$ is the Kronecker delta.}

\section{4. Solution of the inverse problem}

In order to formulate a uniqueness theorem for Inverse Problem 1.1, together with the boundary value problem
$L=L(q(x),a,\alpha,\beta)$ we consider a problem $\tilde L=L(\tilde q(x),a,\alpha,\beta)$ of the same form
but with a different potential $\tilde q(x).$ We agree that if a certain symbol $\gamma$ denotes an object
related to the problem $L,$ then this symbol with tilde $\tilde\gamma$ denotes the corresponding object
related to $\tilde L.$

\medskip
{\bf Theorem 4.1. }{\it In the non-degenerate case: if $\{\lambda_n\}_{n\ge1}=\{\tilde\lambda_n\}_{n\ge1},$ then $q(x)=\tilde q(x)$ a.e. on
$(0,\pi),$ i.e. the specification of the spectrum uniquely determines the potential.

In the degenerate case: let $k>1$ and there exists an operator $K:L_2(0,a)\to L_2(0,a)$ with invertible $I+K,$ such that
\begin{equation}\label{dd:4-1}
q(a-t)=K(q(a+t)), \quad \tilde q(a-t)=K(\tilde q(a+t)), \quad 0<t<a.
\end{equation}
Then the coincidence of the spectra $\{\lambda_n\}_{n\ge1}=\{\tilde\lambda_n\}_{n\ge1}$ also implies $q(x)=\tilde q(x)$ a.e. on $(0,\pi).$
}

\medskip
{\it Proof.} According to Lemma~3.3 the coincidence of the spectra implies
$\Delta_{\alpha,\beta}(\lambda)\equiv\tilde\Delta_{\alpha,\beta}(\lambda),$ which, in turn, by virtue of Theorem~2.1 gives
$W_{\alpha,\beta}(t)=\tilde W_{\alpha,\beta}(t)$ a.e. on $(0,\pi).$ In the non-degenerate case the matrix $A_{\alpha,\beta}$ is invertible.
Thus, according to \eqref{2.7} and invertibility of the operator $R$ we arrive at $q(x)=\tilde q(x)$ a.e. on $(0,\pi).$

In the degenerate case, by virtue of \eqref{2.3.1}, the condition \eqref{dd:4-1} is equivalent to the condition
\begin{equation}\label{dd:4-2}
\left.\begin{array}{l} R_kq(t)=K(R_{k-1}q(t)), \quad R_k\tilde q(t)=K(R_{k-1}\tilde q(t)) \quad \text{for even}\; k,\\[3mm]
R_kq(a-t)=K(R_{k-1}q(a-t)),\quad R_k\tilde q(a-t)=K(R_{k-1}\tilde q(a-t)) \quad \text{for odd}\; k
\end{array}\right\}
\end{equation}
with $t\in(0,a).$ Using the last row in the matrix equality \eqref{sys} along with \eqref{dd:4-2} and invertibility of the operator $I+K$
we obtain
\begin{equation}\label{dd:4-3}
R_{k-1}q(t)=R_{k-1}\tilde q(t), \quad R_kq(t)=R_k\tilde q(t).
\end{equation}
For $j=\overline{1,k-2}$ using $(j+1)$-th row in \eqref{sys} we get the recurrent relations
$$
R_jq(t)=2c(-1)^{\alpha\beta}Q_{j+1}W_{\alpha,\beta}(t)-bcR_{j+2}q(t), \quad j=k-2,k-3,\ldots,1,
$$
which together with \eqref{dd:4-3} and invertibility of $R$ give $q(x)=\tilde q(x)$ a.e. on $(0,\pi).$  $\hfill\Box$

\medskip
{\bf Remark 4.1.} Alternatively to \eqref{1.3}, one can use also other restrictions on $q(x)$ under which the
uniqueness theorem holds in the degenerate case. However, in general, such restrictions depend on the
parameters $\alpha,$ $\beta$ and $k.$ For example, one can use the condition
$$
q(\pi-a+t)=K(q(\pi-a-t)), \quad 0<t<a,
$$
where $I+K$ is invertible for subcases (i) and (iii), while for subcase (ii) the invertibility should be
required for the operator $I-K.$ In particular, this condition includes the limitation to even potentials in
subcases (i) and (iii) or to odd ones in the case (ii), but not vice versa.

\medskip
{\bf Remark 4.2.} If $k\notin{\mathbb N},$ then for the case $\alpha=0,$ $\beta=1$ (corresponding to
non-degenerate subcase (iv) when $k\in{\mathbb N})$ the uniqueness theorem may fail. Indeed, from \eqref{2.6}
and \eqref{allW} it follows that, if, for example, $k=5/2$ and
$$
q(x)=\left\{
\begin{array}{l}
\displaystyle 1,\quad x\in\Big(0,\frac{2\pi}{5}\Big)\cup\Big(\frac{4\pi}{5},\pi\Big),\\[3mm]
\displaystyle -1, \quad x\in\Big(\frac{2\pi}{5},\frac{4\pi}{5}\Big),
\end{array}
\right.
$$
then the problem $L(q(x),2\pi/5,0,1)$ has the same spectrum as the problem $L(0,a,0,1)$ does.

\medskip
For obtaining necessary and sufficient conditions for solvability of Inverse Problem 1.1 we need the following auxiliary result.

\medskip
{\bf Lemma 4.1. }{\it Fix $\alpha,\beta\in\{0,1\}.$ Let arbitrary complex numbers $\lambda_n,$ $n\ge1,$ of
the form \eqref{dd:5} be given. Then there exists a function $W_{\alpha,\beta}(t)\in L_2(0,\pi)$ such that
the function $\Delta_{\alpha,\beta}(\lambda)$ determined by \eqref{EQ} has the form \eqref{2.5} or
\eqref{2.6} depending on $\alpha=\beta$ or $\alpha\ne\beta,$ respectively.

Moreover, for the combinations of $\alpha,$ $\beta$ and $k$ satisfying one of the groups of conditions in
\eqref{deg}, if additionally the corresponding condition in \eqref{dd:6}--\eqref{BBV:9} of the degeneration
of the numbers $\lambda_n$ is fulfilled, then the involved function $W_{\alpha,\beta}(t)$ satisfies one of
the conditions \eqref{dd:3}--\eqref{nn:3}, respectively.}

\medskip
{\it Proof.} Let $\alpha=\beta=0$ and the sequence of complex numbers $\{\lambda_n\}_{n\ge1}$ of the form
\eqref{dd:5} be given. By the standard approach (see, e.g., \cite{but1}) it can be proven that the function
$\Delta_{0,0}(\lambda)$ constructed via \eqref{EQ} has the form \eqref{2.5} with $\alpha=0.$ It remains to
prove \eqref{dd:3}.

By virtue of the entireness of the function $\Delta_{0,0}(\lambda)$ we have
\begin{equation}\label{dd:11-1}
\int_0^{\pi} W_{0,0}(t)\,dt=0.
\end{equation}
Further, assuming \eqref{dd:6} or, in other words, $\Delta_{0,0}((kn)^2) = 0$ for $n \in \mathbb N$ and
substituting $\lambda_{kn} = (kn)^2$ into \eqref{2.5}, we get
\begin{equation}\label{dd:11}
\int_0^{\pi} W_{0,0}(t) \cos k n t \, dt = 0, \quad n \in \mathbb N.
\end{equation}
Expanding $W_{0,0}(t)$ into the Fourier series \eqref{dd:7}, we obtain $a_0 = 0$ and $a_{k n} = 0$, $n \in
\mathbb N$, by virtue of \eqref{dd:11-1} and \eqref{dd:11}, respectively. Thus, we obtain
$$
\sum_{n=0}^\infty a_nT_n\cos nt=0,
$$
where $T_n$ is determined in \eqref{dd:8-1}, and then according to \eqref{dd:8} we arrive at \eqref{dd:3}.
Degenerate subcases (ii) and (iii) are treated analogously. $\hfill\Box$

\medskip
The following theorem means that the properties of the spectrum proven in Theorem~3.1, actually being necessary, are also sufficient
conditions for the solvability of Inverse Problem 1.1. In other words, Theorem~3.1 gives a full characterization of the spectrum of the
boundary value problem $L$ both in the degenerate and non-degenerate cases.

\medskip
{\bf Theorem 4.2. }{\it (I) Non-degenerate case. Let $\alpha,$ $\beta$ and $k$ satisfy one of the groups of conditions in \eqref{non-deg}.
Then for an arbitrary sequence of complex numbers $\{\lambda_n\}_{n\ge1}$ of the form \eqref{dd:5} there exists a function $q(x)\in
L_2(0,\pi)$ such that $\{\lambda_n\}_{n\ge1}$ is the spectrum of the boundary value problem $L(q(x),\pi/k,\alpha,\beta).$

(II) Degenerate case. Let $\alpha,$ $\beta$ and $k$ satisfy one of the groups of conditions in \eqref{deg}. Then for an arbitrary sequence
of complex numbers $\{\lambda_n\}_{n\ge1}$ of the form \eqref{dd:5} satisfying the corresponding degeneration condition in
\eqref{dd:6}--\eqref{BBV:9} there exists a function $q(x)\in L_2(0,\pi)$ (not unique) such that $\{\lambda_n\}_{n\ge1}$ is the spectrum of
the boundary value problem $L(q(x),\pi/k,\alpha,\beta).$}

\medskip
{\it Proof.} Fix $\alpha,\beta\in\{0,1\}$ and $k\in{\mathbb N}.$ Using the given sequence $\{\lambda_n\}_{n\ge 1}$ of the form \eqref{dd:5}
construct the function $\Delta_{\alpha,\beta}(\lambda)$ by formula  \eqref{EQ}. According to Lemma 4.1, $\Delta_{\alpha,\beta}(\lambda)$
has the form \eqref{2.5} or \eqref{2.6} with a certain function $W_{\alpha,\beta}(t)\in L_2(0,\pi).$

(I) Under any assumption in \eqref{non-deg}, according to Lemma~2.1, the main equation \eqref{2.7} has a unique solution $q(x)\in
L_2(0,\pi).$ Consider the boundary value problem $L:=L(q(x),\pi/k,\alpha,\beta)$ with this $q(x).$ It can be easily seen that the spectrum
of $L$ coincides with $\{\lambda_n\}_{n\ge 1}.$

(II) For $k=1$ the assertion is obvious. Let $k>1$ and one of the conditions in \eqref{deg} along with the corresponding degeneration
condition in \eqref{dd:6}--\eqref{BBV:9} be fulfilled. Then, by virtue of Lemma~4.1, the function $W_{\alpha,\beta}(t)$ satisfies the
corresponding condition in \eqref{dd:3}--\eqref{nn:3}. Let us show the solvability of the main equation \eqref{2.7}, which, in turn, is
equivalent to \eqref{sys}. For briefness we rewrite \eqref{sys} in the form
\begin{equation}\label{last-1}
Y=A_{\alpha,\beta}X,
\end{equation}
where $Y=(y_1,\ldots,y_k)^T,$
$y_j=2(-1)^{\alpha\beta}Q_jW_{\alpha,\beta}(t)$ and $X=(x_1,\ldots,x_k)^T$ is unknown. Consider for
definiteness case (i) in \eqref{deg}, i.e. $\alpha=\beta=0.$ Other cases (ii) and (iii) in \eqref{deg} can be
treated similarly. Then the corresponding condition \eqref{dd:3} is equivalent to
\begin{equation}\label{rank}
\sum_{j=1}^ky_j=0.
\end{equation}
According to the proof of Lemma~2.1, $\rank A_{0,0}=k-1.$ Since the sum of rows of the matrix $A_{0,0}$ is the zero-row, the relation
\eqref{rank} gives $\rank [Y,A_{0,0}]=k-1,$ i.e. the rank of the extended matrix $[Y,A_{0,0}]$ equals to the rank of $A_{0,0}.$ Thus, the
system $Y=A_{0,0}X$ has a solution (not unique). Since $A_{0,0}$ is a numerical matrix, the elements of $X$ and $Y$ can belong only to one
and the same linear space (in our case $L_2(0,a)).$ Hence, $X\in(L_2(0,a))^k.$ Put $q(x)=R^{-1}X$ and consider the problem
$L(q(x),\pi/k,\alpha,\beta)$ with this $q(x).$ Obviously, $\{\lambda_n\}_{n\ge 1}$ is its spectrum. $\hfill\Box$

\medskip
The proof of Theorem 4.2 is constructive and gives algorithms for solving Inverse Problem~1.1. The following algorithm allows one to
construct the solution of the inverse problem in the non-degenerate case.

\medskip
{\bf Algorithm 4.1.} Let the spectrum $\{\lambda_n\}_{n\ge 1}$ of a boundary value problem $L(q(x),a,\alpha,\beta)$ in the non-degenerate
case be given.

1. Construct $\Delta_{\alpha,\beta}(\lambda)$ by formula \eqref{EQ}.

2. Find $W_{\alpha,\beta}(t)$ in the corresponding representation \eqref{2.5} or \eqref{2.6} inverting the Fourier transform.

3. Construct $q(x)$ as a solution of the main equation \eqref{2.7}.

\medskip
Lemma~3.2 along with formula \eqref{allW} give the following algorithm, which implements the third step in Algorithm~4.1.

\medskip
{\bf Algorithm 4.2.} Let the function $W_{\alpha,\beta}(x)$ in the non-degenerate case be given.

1. Calculate $q(x)$ on $(0,a)$ by the corresponding formula in Lemma~3.2.

2. Construct $q(x)$ on $(a,2a)$ by the formula
$$
q(x)=2(-1)^\alpha W_{\alpha,\beta}(\pi+a-x)-q(2a-x), \quad x\in(a,2a).
$$

3. For $j=\overline{2,k-1}$ repeat the following step. Let the function $q(x)$ be already calculated on the interval $(0,ja).$ Then find
$q(x)$ on $(ja,(j+1)a)$ by the formula
$$
q(x)=(-1)^\alpha\Big(2W_{\alpha,\beta}(\pi+a-x)+q(x-2a)\Big), \quad x\in(ja,(j+1)a).
$$

\medskip
Let us now give an algorithm for solving the inverse problem in the degenerate case. For definiteness we assume that $k>1$ and the
potential $q(x)$ satisfies the condition \eqref{1.3} with some known operator $K.$

\medskip
{\bf Algorithm 4.3.} Let the spectrum $\{\lambda_n\}_{n\ge 1}$ of a boundary value problem $L(q(x),a,\alpha,\beta)$ in the degenerate case
along with the operator $K$ in \eqref{1.3} be given.

1. Construct $\Delta_{\alpha,\beta}(\lambda)$ by formula \eqref{EQ}.

2. Find $W_{\alpha,\beta}(t)$ in the corresponding representation \eqref{2.5} or \eqref{2.6} inverting the Fourier transform.

3. Calculate $q(x)$ on the interval $(a,2a)$ by the formula
$$
q(a+x)=(I+K)^{-1}(-2W_{\alpha,\beta}(\pi-x)), \quad x\in(0,a).
$$

4. Calculate $q(x)$ on the interval $(0,a)$ by the formula
$$
q(a-x)=-2W_{\alpha,\beta}(\pi-x)-q(a+x), \quad x\in(0,a).
$$

5. For $j=\overline{2,k-1}$ repeat the following step. Let the function $q(x)$ be already calculated on the interval $(0,ja).$ Then find
$q(x)$ on $(ja,(j+1)a)$ by the formula
$$
q(x)=-2W_{\alpha,\beta}(\pi+a-x)+(-1)^\alpha q(x-2a), \quad x\in(ja,(j+1)a).
$$

\medskip
In the degenerate case Algorithm 4.3 allows one to describe the set of all iso-spectral potentials $q(x),$
i.e. of those  for which the corresponding problems $L(q(x),a,\alpha,\beta)$ have one and the same spectrum
$\{\lambda_n\}_{n\ge 1}.$ For this purpose on the third step of the algorithm one should use a constant
operator $K,$ i.e. when there exists a function $p(x)\in L_2(0,a)$ such that
\begin{equation}\label{last}
K(f(x))=p(x)
\end{equation}
for all $f(x)\in L_2(0,a).$ Indeed, the following theorem holds.

\medskip
{\bf Theorem 4.3. }{\it If the function $p(x)$ in \eqref{last} varies through $L_2(0,a),$ then the functions
$q(x)$ constructed by Algorithm~4.3 form the set of all iso-spectral potentials for the given spectrum
$\{\lambda_n\}_{n\ge 1}.$ }

\medskip
{\it Proof.} It is clear that for the operator $K$ of the form \eqref{last} and any $p(x)\in L_2(0,a),$ Algorithm~4.3 gives iso-spectral
potentials $q(x)$ with $q(x)=p(a-x)$ a.e. on $(0,a).$ On the other hand, from Theorem~4.1 it follows that no other iso-spectral potentials
exist. $\hfill\Box$

\medskip
{\bf Remark 4.3.} For describing the set of iso-spectral potentials one can arbitrarily specify $q(x)$ on the fixed interval $((j-1)a,ja)$
not only for $j=1$ but also for any fixed $j=\overline{2,k}.$ Indeed, one can easily check that in the matrix $A_{\alpha,\beta}$ for any
$j=\overline{1,k}$ there exists a basis minor that does not include elements of $j$-th column. Thus, in the degenerate linear system
\eqref{last-1} the variable $x_j$ can be considered as a free one. Consequently, once being solvable this system remains to be so for any
preassigned value of $x_j.$

\bigskip
{\bf Acknowledgement.} This research was supported in part by RFBR (Grants 15-01-04864, 16-01-00015 and 17-51-53180) and by the Ministry of
Education and Science of RF (Grant 1.1660.2017/4.6). The co-author N.P.~Bondarenko was also supported by Russian Federation President Grant
MK-686.2017.1.

\end{document}